\documentclass[12pt]{article}
\usepackage{amsfonts}
\usepackage{mathrsfs}
\usepackage{amsmath,amssymb}

\baselineskip 5.5pt \lineskip 5.5pt \numberwithin{equation}{section}

\def\QED{\hfill$\Box$\par}
\def\cl{\centerline}
\def\vs{\vspace*}
\def\ni{\noindent}
\def\C{\mathbb{C}}
\def\Z{\mathbb{Z}}
\def\la{\lambda}

\def\a{\alpha}

\def\bo{{\bf 1}}
\def\rar{\rightarrow}
\def\ww{\mathcal{W}}
\def\uu{\mathcal{U}}
\def\pp{\mathcal{P}}
\def\pq{\mathcal{P}_{+}}
\def\p-{\mathcal{P}_{\leq 0}}

\begin{document}
\cl{{\bf\large WHITTAKER MODULES FOR THE $W$-ALGEBRA $W(2,2)$}}
\footnote {\!\!\!The second author is supported by the Postdoctoral
Sience Foundation of China (Grant No. 20080440720) and by NSF
 of China (Grant No. 10671027)} \vs{6pt}

\cl{Bin Wang, \,\,  Junbo Li}

\begin{abstract}{\footnotesize In this paper, Whittaker modules for the $W$-algebra
$W(2,2)$ are studied. We obtain analogues to several results from
the classical setting and the case of the Virasoro algebra,
including a classification of simple Whittaker modules by central
characters and composition series for general Whittaker modules.}
\end{abstract}


\noindent{\it{MR(2000) Subject Classification}: 17B10, 17B65,
17B68.}\vs{10pt}

\cl{\bf\S1.Introduction}

In this paper we investigate Whittaker modules for the $W$-algebra
$W(2,2)$.  Whittaker modules were first discovered for
$\mathfrak{sl}_2{(\C )}$ by Arnal and Pinzcon in \cite{AP}. Block
showed, in \cite{B} that the simple modules for $\mathfrak{sl}_2(\C
)$ consist of highest (lowest) weight modules, Whittaker modules and
a third family obtained by localization. This illustrates the
prominent role played by Whittaker modules. The algebra $W(2,2)$ is
related to vertex operator algebras. In \cite{ZD}, $W(2,2)$ plays an
important role in the classification of simple vertex operator
algebras generated by two weight vectors.

Kostant defined Whittaker modules for an arbitrary
finite-dimensional complex semisimple Lie algebra $\mathfrak{g}$ in
\cite{K}, and showed that these modules, up to isomorphism, are in
bijective correspondence with ideals of the center
$Z(\mathfrak{g})$. In particular, irreducible Whittaker modules
correspond to maximal ideals of $Z(\mathfrak{g})$. In the quantum
setting, Whittaker modules have been studied by Sevoystanov for
$\uu_h(\mathfrak{g})$ in \cite{S} and by M. Ondrus for
$\uu_q(\mathfrak{sl}_2)$ in \cite{O}. Recently Whittaker modules
have also been studied by M. Ondrus and E. Wiesner for the Virasoro
algebra in \cite{OW}, X. Zhang and S. Tan for the
Schr\"{o}dinger-Virasoro algebra in \cite{ZT}, K. Christodoulopoulou
for Heisenberg algebras in \cite{C}, and by G. Benkart and M. Ondrus
for generalized Weyl algebras in \cite{BO}.

The algebra $W(2,2)$ is an infinite-dimensional Lie algebra.
Recently, there appeared some papers on the structures and
representations for this algebra. In \cite{LGZ}, \cite{LZ} and
\cite{ZD}, its irreducible weight modules, indecomposable modules
and Verma modules were respectively studied. Its derivations,
central extensions and automorphisms were determined in \cite{GJP}.
The Lie bialgebra and left symmetric algebra structures on the
algebra were determined in \cite{LS} and \cite{CL}.

Our main result is Theorem 4.2 in which we get a classification of
irreducible Whittaker modules by central characters. This theorem
follows easily from two results, one is Proposition 3.1 by which the
Whittaker vectors in a universal Whttaker module are determined, and
the other is Lemma 3.3 which says that every submodule of a
universal Whittaker module contains a Whittaker vector. We use the
concrete nature of the algebra $W(2,2)$ to obtain 3.1, but our proof
for 3.3, different from \cite{OW}, is more general and could be
possibly extended to infinitely generated algebras. The paper is
organized as follows.  In section 2, we define Whittaker vectors and
Whittaker modules for $W(2,2)$, and also construct a universal
Whittaker module. Then the Whittaker vectors in a universal
Whittaker module are examined in section 3 and the irreducible
Whittaker modules are classified in section 4. In the last section
we describe a decomposition of an arbitrary Whittaker module and
characterize its submodules.

\vs{18pt}

\cl{\bf\S2. Preliminaries}

\ni 2.1. The algebra $W(2,2)$, denoted by $\ww$, is an infinitely
dimensional Lie algebra with a $\C$-basis
$\{\,L_n,\,W_n,\,z\,|\,n\in \Z\,\}$ and the following Lie brackets:
\begin{eqnarray*}\label{BLie}\begin{array}{lll}
&&[L_n,L_m\,]=(m-n)L_{m+n}+\frac{n^3-n}{12}\delta_{m+n,0}z,\\[6pt]
&&[L_n,W_m]=(m-n)W_{m+n}+\frac{n^3-n}{12}\delta_{m+n,0}z.\\
&&[W_n,W_m]=[z,W_m]=[z,L_m]=0\end{array}
\end{eqnarray*}
Let $S(z)$represent the symmetric algebra generated by $z$, that is,
plolynomials in $z$. Then $S(z)$ is obviously contained in the
center of the universal enveloping algebra $\uu (\ww )$.

Set ${\ww}_n=Span_{\C}\{L_n, W_n\},n\neq 0,{\ww}_0=Span_{\C}\{L_0,
W_0,z\}$. Then $\ww =\bigoplus_{n\in \Z}{\ww}_n$ and it is easy to
check that $[{\ww}_n, {\ww}_m]\subset {\ww}_{n+m}$, i.e. $\ww$ is a
$\Z$-graded Lie algebra. Furthermore, set
${\ww}_{-}=\bigoplus_{n<0}{\ww}_n, {\ww}_{+}=\bigoplus_{n>0}{\ww}_n,
{\ww}_{\leq 0}=\bigoplus_{n\leq 0}{\ww}_n$, then clearly ${\ww}_{-},
{\ww}_0, {\ww}_{+}$ and ${\ww}_{\leq 0}$ are subalgebras of $\ww$.

\ni 2.2. {\bf{Partitions}}.

\ni 2.2.1. We define a partition to be a non-decreasing sequence of
integers
$\mu=(\mu_1,\mu_2,\cdots,\mu_r),\mu_1\leqslant\mu_2\leqslant\cdots\leqslant\mu_r.$
Denote by $\pp$ the set of all partitions. For $\lambda
=(\la_1,\cdots ,\la_r) \in \pp$, we define the length of $\la$ to be
$r$, denoted by $\ell (\la)$; if all $\la_i >0, 1\leq i \leq r$,
then call $\la$ a positive partition, and if all $\la_i \leq 0$ then
call $\la$ a non-positive partition. Denote by $\p-$ the set of all
non-positive partitions and by $\pq$ the set of all positive
partitions.

For $\la=(\la_1,\cdots ,\la_r) \in \pp ,k\in \Z$, let $\la (k)$
denote the number of times $k$ appears in the partition. Clearly the
values $\la (k), k\in \Z$ completely determine the partition $\la$.
So, we sometimes write $\la$ in an alternative form, $\la =(\cdots
(-1)^{\la (-1)},0^{\la (0)},1^{\la (1)},\cdots )$. Note that $\la
(k)=0 $ when $|k|$ is sufficiently large. Define elements $L_{\la},
W_{\la}\in \uu{(\ww )}$ by
\begin{eqnarray*}
&&L_{\la}=L_{\la_1}L_{\la_2}\cdots L_{\la_r}=\cdots
L_{-1}^{\la{(-1)}}L_{0}^{\la{(0)}}L_{1}^{\la{(1)}}\cdots\\
&&W_{\la}=W_{\la_1}W_{\la_2}\cdots W_{\la_r}=\cdots
W_{-1}^{\la{(-1)}}W_{0}^{\la{(0)}}W_{1}^{\la{(1)}}\cdots
\end{eqnarray*}
Set $\bar{0}=(\cdots (-1)^0,0^0,1^0,\cdots )$ and
$L_{\bar{0}}=W_{\bar{0}}=1 \in \uu{(\ww )}$. We will consider
$\bar{0}$ to be an element of $\p-$ but not of $\pq$. By PBW
theorem, we know that $\{L_{\la}W_{\la'}L_{\mu}W_{\mu'}|\la ,\la'
\in \p- ,\mu ,\mu' \in \pq \}$ form a basis of $\uu (\ww )$ over
$S(z)$.

\ni 2.2.2. $\uu{(\ww )}$ naturally inherits a grading from the one
of $\ww$. Namely, set $\uu{(\ww )}_m =Span_{\C}\{x_1\cdots x_k |
x_i\in \ww_{n_i}, 1\leq i\leq k, \sum_{i=1}^{k}n_i =m \}$, and then
$\uu{(\ww )}=\bigoplus_{m\in \Z}\uu{(\ww )}_m$ is a $\Z$-graded
algebra, i.e. $\uu{(\ww )}_m \uu{(\ww )}_n \subset \uu{(\ww
)}_{m+n}$. Similarly, $\uu (\ww_{+} )$ (resp. $\uu (\ww_{\leq 0})$)
inherit $\Z_+$-grading (resp.$\Z_{\leq 0}$-grading) from $\ww_+$,
(resp. $\ww_{\leq 0})$. If $x \in \uu{(\ww )}_m$, then we say $x$ is
a homogeneous element of degree $m$. If set $|\la |=\la_1 + \la_2
+\cdots +\la_{{\ell}(\la)}$, $L_{\la}$ and $W_{\la}$ are homogeneous
elements of degree $|\la |$. If $x (\neq 0 )$ is not homogeneous but
a sum of finitely many nonzero homogeneous elements, then denote by
$mindeg(x)$ the minimum degree of its homogeneous components.
Moreover, let us, for convenience, call any product of elements
$L_m^t ,W_n^k ,z^q ,(m,n \in \Z ,t,q,k \geq 0)$, in $\uu{(\ww )}$ a
monomial, of height (resp. height w.r.t. $L$ ) equal to the sum of
the various $t$'s and $k$'s (resp. $t$'s ) occurring. If $ x \in
\uu{(\ww )}$ is a sum of monomials of height (resp. height w.r.t. L
) $\leq l$, we write $ht(x)\leq l$ (resp. $ht_1(x) \leq l)$. Then we
have, by PBW theorem,

\ni {\bf Lemma}. \ \ {\it For $m,n \in \Z$ ,let $A_m$ stands for
either $L_m$ or $W_m$ (resp. $L_m$ ), then $A_m^t A_n^k$ is a linear
combination of $A_n^k A_m^t$ along with other monomials of height
(resp. height w.r.t. $L$ ) $< t+k$.} \QED

\ni 2.2.3. We need some more notation. For $\la =(\la_1,\la_2,\cdots
\la_r)\in \pp ,1\leq i\leq r, 0\leq j < r, $write
\begin{eqnarray*}
&&\la \{i\}=(\la_1 ,\cdots , \la_i ), \,\, \la \{0\}=\bar{0}, \\
&&\la [j]=(\la_{j+1},\cdots ,\la_r ), \,\, \la [r]=\bar{0},\\
&&\la <i>=(\la_1,\cdots ,\la_{i-1},\hat{\la_i},\la_{i+1},\cdots ).
\end{eqnarray*}

\ni{\bf Lemma A} \ \  {\it $ht_1 ([W_m, L_{\la}]) < \ell (\la),
\forall m \in \Z, \la \in \pp$.}

\ni {\bf Proof}\ \
$[W_m,L_{\la}]=\mbox{$\sum\limits_{i=1}^{{\ell}(\la)}$}L_{\la
\{i-1\}}[W_m, L_{\la_i}]L_{\la [i]}$. Note there is no $L_k$
appearing in $[W_m, L_{\la_i}]$ and hence the lemma follows.\QED

\ni {\bf Lemma B}\ \ {\it  Write $\uu'$ for $\uu ({\ww_{+}})$,
$\uu''$ for $\uu (\ww_{\leq 0})$. Let $0\neq x \in \uu'_n, 0\neq y
\in \uu''_m $ with $n
> 0, m \leq 0$ and $s = max\{m+n, 0 \}$, then $[x, y] =
\sum_{k=s}^{n} u_k $ with $u_k = \sum\limits_{i}y_{k,i}x_{k,i}$
where $x_{k,i} \in \uu'_k, y_{k, i} \in \uu'' _{n+m-k}$ and
$ht(y_{k,i}) < ht(y)$ if $k = n, y_{k,i} \neq 0$. Further, if assume
$x \in (\ww_{+})_n$, then $x_{k,i} \in (\ww_{+})_k$ and $ht(y_{k,i})
< ht(y)$ whenever $k > 0, y_{k,i} \neq 0$. }

\ni {\bf Proof}\ \ By direct checking.\QED

\ni 2.3 {\bf{Whittaker Module}}. Observe that $\ww_+ $ is a
subalgebra of $\ww$ generated by $L_1,L_2,W_1,W_2$ As in \cite{OW},
we define Whittaker modules for $\ww$ as what follows.

\ni 2.3.1. {\bf{Definition}}. One says a Lie algebra homomorphism
$\varphi:\ww_+ \rar \C$ non-singular, if
$\varphi(L_i)\varphi(W_i)\neq0$, $i=1,2$. For a $\ww$-module $V$, a
vector $v\in V$ is said to be a Whittaker vector if $xv=\varphi(x)v$
for all $x\in \ww_+$. Furthermore, if $v$ generates $V$, we call $V$
a Whittaker module of type $\varphi$ ($\varphi$ is required to be
non-singular in this case) and $v$ a cyclic Whittaker vector of $V$.

\ni 2.3.2. For a given non-singular Lie algebra homomorphism
$\varphi:\ww_+\rar\C$, define $\C_\varphi$ to be the one-dimensional
$\ww_+$-module given by the action $x\a= \varphi(x) \a$ for all
$x\in \ww_+$ and $\a\in\C$. Then the induced $\ww$-module
\begin{eqnarray*}
&&M_\varphi=\uu(\ww)\otimes_{\uu(\ww_+)}\C_\varphi,
\end{eqnarray*}
is a Whittaker module of type $\varphi$ with a cyclic Whittaker
vector $w=\bo\otimes 1$.  By PBW theorem, it's easy to see that
$\{z^kL_{\la}W_{\mu}w\,|\,\la,\mu\in\p-,k\geq 0\}$ is a basis of
$M_\varphi$ as a vector space over $\C$.

Besides, for any  $\xi\in\C$, obviously $(z-\xi)M_\varphi$ is a
submodule of $M_\varphi$. Define
$L_{\varphi,\xi}=M_\varphi/(z-\xi)M_\varphi$ and denote the
canonical homomorphism by $p_\xi$. Then $L_{\varphi,\xi}$ is a
Whittaker module for $\ww$. The following lemma makes $M_{\varphi}$
become a universal Whittaker module.

\ni {\bf Lemma} {\it Fix $\varphi$ and $M_{\varphi}$ as above. Let V
be a Whittaker module of type $\varphi$ generated by a Whittaker
vector $w_v$. The there is a unique map $\phi:M_{\varphi}\rightarrow
V$ taking $w=1\otimes 1$ to $w_v$.}

\ni{\bf Proof}.\ \ Uniqueness is obvious. For $u\in \uu(\ww)$, one
can write, by PBW,
$$u=\sum\limits_\alpha b_\alpha n_\alpha,b_\alpha\in \uu(\ww_{\leq 0}),\,\, n_\alpha\in U(W_+).$$
If $uw=0$ then $uw=\sum\limits_\alpha b_\alpha \psi(n_\alpha)w=0$.
Therefore, $\sum\limits_\alpha b_\alpha \psi(n_\alpha)=0$. Now it's
easy to see that the map $\phi : M_{\varphi} \rar V$, defined by
$\phi(uw) = uw_v$, is well defined.\QED

\ni 2.3.3.  Write $M$ for $M_{\varphi}$, $\uu'$ for $\uu (\ww_{+})$,
and $\uu''$ for $\uu (\ww_{\leq 0})$. For $k \leq 0, $ define $M_k =
\{ xw \,|\, x \in {\uu''_k} \}$ and then clearly $M =M_0 \oplus
M_{-1} \oplus M_{-2} \oplus \cdots$. We say that a nonzero
homogeneous vector $v$ in $M$ is of degree k if $v \in M_k$. If $v
(\neq 0 )$ is not homogeneous but a sum of finitely many nonzero
homogeneous vectors, then define $mindeg(v)$ to be the minimum
degree of its homogeneous components. Meanwhile, for any nonzero
vector $v \in M$, let $d = mindeg(v)$ and then there uniquely exist
$v_i \in \uu''_i, 0 \geq i \geq d$ such that $v =
\sum\limits_{i=d}^{0} v_{i}w$ with $v_{d} \neq 0$. Then define $\ell
(v) = ht( v_{d})$ and $\ell'(v) = ht_1 (v_{d})$.

\vs{18pt}

\cl{\bf\S3. Whittaker Vectors in $M_{\varphi}$ and $L_{\varphi
,\xi}$}

In this section, we characterize the Whittaker vectors in
$M_{\varphi}$ and $L_{\varphi ,\xi}$. Fix a nonsingular Lie algebra
homomorphism $\varphi : \ww_+ \rar \C$, and let $w=\bo\otimes 1 \in
M_{\varphi}$.

\ni 3.1. {\bf Proposition} \ \ {\it If $w'\in M_\varphi$ is a
Whittaker vector (of type $\varphi$), then $w'=p(z)w$ for some
$p(z)\in S(z).$}

\ni{\bf Proof}. \ \ Write $w'=\sum\limits_{\lambda,\mu\in
\p-}p_{\lambda,\mu}(z)L_{\lambda}W_{\mu}w$. Clearly it's enough to
show that $p_{\la,\mu}=0$ unless $\la =\bar{0}$ and $\mu =\bar{0}$.

a). Suppose there, at least, exists a $\lambda\neq \bar{0}$ such
that $p_{\lambda,\mu}(z)\neq 0$ for some $\mu$. Now let $A=\{\la
\,|\, p_{\la,\mu}\neq 0, \la\in \p- \},\, l=max\{\ell(\la) \,|\,
\la\in A\},$ and $B=\{\la | \ell(\la )=l,\la \in A\}$. Put
$m_0=min\{\la_1 \,|\, \la\in B\}$ and $B'=\{\la \in B \,|\,
\la_1=m_0 \}$. Note that if $\la \in B, \la_i =m_0$ for some $i$,
then $\la_1=\cdots =\la_i=m_0$ and hence $\la \in B', \,
\la<i>=\la<1>$.

Let $m=2-m_0 (\geq 2)$ and then $\forall \la \in B, m+\la_i >2$
unless $\la_i =m_0 $. Consider
\begin{eqnarray*}
(W_m -\varphi(W_m))w'&=&\mbox{$\sum\limits_{\la,\mu}$}p_{\la,\mu}(z)
[W_{m},L_{\la}W_{\mu}]w\\
&=&\mbox{$\sum\limits_{\la,\mu}$}p_{\la,\mu}(z)
[W_{m},L_{\la}]W_{\mu}w=D_1 + D_2.
\end{eqnarray*}
where $D_1 =\mbox{$\sum\limits_{\mu,\la \in B}$}p_{\la,\mu}(z)
[W_{m},L_{\la}]W_{\mu}w$, and $D_2 =\mbox{$\sum\limits_{\mu,\la
\notin B}$}p_{\la,\mu}(z) [W_{m},L_{\la}]W_{\mu}w$. Then, clearly
$$\ell'(p_{\la,\mu}(z) [W_{m},L_{\la}]W_{\mu}w) \leq l-2,$$ for any $
\la \notin B$ and hence $\ell'{(D_2)}\leq l-2$, by 2.2.3 Lemma A.
Meanwhile,
\begin{eqnarray*}
D_1 &=&\mbox{$\sum\limits_{\mu,\la \in
B}$}p_{\la,\mu}(z)\mbox{$\sum\limits_{i=1}^{\ell(\la
)}$}L_{\la\{i-1\}} [W_{m},L_{\la_i}] L_{\la[i]}W_{\mu}w\\
&=&\mbox{$\sum\limits_{\mu,\la \in
B}$}p_{\la,\mu}(z)\mbox{$\sum\limits_{i=1}^{\ell(\la
)}$}L_{\la\{i-1\}}(\la_i - m)W_{m+\la_i}L_{\la[i]}W_{\mu}w\\
 &=& D_1' +D_1'',
\end{eqnarray*}
where \begin{eqnarray*} &&D_1' =\mbox{$\sum\limits_{\mu,\la \in
B}$}p_{\la,\mu}(z)\mbox{$\sum\limits_{i=1}^{\ell(\la
)}$}L_{\la<i>}(\la_i - m)W_{m+\la_i} W_{\mu}w,\\
&&D_1'' =\mbox{$\sum\limits_{\mu,\la \in
B}$}p_{\la,\mu}(z)\mbox{$\sum\limits_{i=1}^{\ell(\la
)}$}L_{\la\{i-1\}}(\la_i - m)[W_{m+\la_i}, L_{\la[i]}] W_{\mu}w.
\end{eqnarray*}
Then, we have $\ell'(D_1'') \leq l-2$, by lemma A. But, since
$W_{m+\la_i}w =0$ when $m+\la_i >2$,
\begin{eqnarray*} D_1' &=&\mbox{$\sum\limits_{\mu,\la \in
B'}$}p_{\la,\mu}(z)\mbox{$\sum\limits_{i=1}^{\la(m_0)
}$}(2m_0-2)L_{\la<i>}W_2 W_{\mu}w\\
&=&\mbox{$\sum\limits_{\mu,\la \in B'}$}(2m_0 - 2)\varphi (W_2)\la
(m_0)p_{\la,\mu}(z)L_{\la <1>}W_{\mu}w.
\end{eqnarray*}
So, $D_1'$ equals to a linear combination of some vectors $v$'s in
$M_{\varphi}$ with $ht_1(v) = l-1$ which are linearly independent
from each other. Therefore $D_1'$ and $D_1'' + D_2$ are linearly
independent. Thus, $(W_m -\varphi(W_m))w'=D_1'+ D_1'' + D_2\neq 0$
which contradicts with the hypothesis that $w'$ is a Whittaker
vector.

b). Suppose $\exists \mu\neq \bar{0}, p_{\bar{0},\mu}(z)\neq 0$, and
$p_{\la,\mu}(z)=0, \forall \la \neq \bar{0}$. Then we write $w'=
\mbox{$\sum\limits_{\mu \in \p-}$}p_{\bar{0},\mu}(z) W_{\mu}w$. Let
$C=\{\mu \,|\, p_{\bar{0},\mu}(z) \neq 0 \}, n_0 = min\{\mu_1 | \mu
\in C \}$, and $ C'= \{\mu \,|\, \mu_1 = n_0, \mu \in C \}$. And if,
$\forall \mu \in C, \mu_i = n_0$ for some $i$, then clearly $\mu_1
=\cdots =\mu_i =n_0$ and $\mu \in C'$. Now let $m= 2-n_0 $ and then
$m+\mu_i >2 $ unless $\mu_i =n_0$.

Consider
\begin{eqnarray*}
(L_m-\varphi(L_m))w' &=&\mbox{$\sum\limits_{\mu}$}p_{\bar{0},\mu}(z)
[L_{m},W_{\mu}]w\\
&=&\mbox{$\sum\limits_{\mu}$}p_{\bar{0},\mu}(z)\sum_{i=1}^{\ell
(\mu)}W_{\{i-1\}}[L_{m},W_{\mu_i}]W_{\mu [i]}w \\
&=&\mbox{$\sum\limits_{\mu}$}p_{\bar{0},\mu}(z)\sum_{i=1}^{\ell
(\mu)}(m-\mu_i )W_{\mu <i>}W_{m+\mu_i}.
\end{eqnarray*}
Since $W_{m+\mu_i}w= 0$ if $m+\mu_i > 2$, we have
$$ (L_m - \varphi (L_m ))w'=\mbox{$\sum\limits_{\mu \in C'}$}p_{\bar{0},\mu}(z)\mu (n_0)(m-n_0 )\varphi{(W_2)}W_{\mu
<1>}w. $$ But this is a sum of terms that are linearly independent
from each other and hence $(L_m - \varphi (L_m ))w' \neq 0$ which
contradicts with the hypothesis that $w'$ is a Whittaker vector.
\QED

\ni 3.2. {\bf Proposition} \ \ {\it Let $w=\bo \otimes1\in
M_{\varphi}$ and $\overline{w}=\overline{\bo \otimes1}\in
L_{\varphi,\xi}$. If $w'\in L_{\varphi,\xi}$ is a Whittaker vector
then $w'=c\overline{w}$ for some $c \in \C$.}

\ni{\bf Proof} \ \ $L_{\varphi,\xi}$ has a basis in the form
$\{L_{\la}W_{\mu}\overline{w}\ | \la, \mu\in \p-\}$. In fact,
clearly it is enough to show this set is linearly independent.
Suppose there are $a_{\la,\mu}\in \C$, with at most finitely many
$a_{\la,\mu}\neq 0$, such that
$$0 = \mbox{$\sum\limits_{\la,\mu \in \p-}$}a_{\la,
\mu}L_{\la}W_{\mu}\overline{w} = \overline{\mbox{$\sum\limits_{\a,
\mu}$}a_{\la, \mu}L_{\la}W_{\mu}w}$$in $L_{\varphi, \xi }$. So
$$\mbox{$\sum\limits_{\mu}$}a_{\la, \mu}L_{\la}W_{\mu}w = (z-\xi
)\mbox{$\sum\limits_{\la,\mu}\sum\limits_{i=1}^{k}$}b_{\la,
\mu,i}z^iL_{\la}W_{\mu}w$$ for some $k > 0$ and $b_{\la,\mu,i} \in
\C$. The right hand side of this expression can be rewritten as a
linear combination of $\C$-basis vectors $z^i L_{\la}W_{\mu}w$ in
$M_{\varphi, \xi}$. Then comparing all the coefficients of the two
sides of the above equation, we can deduce that $a_{\la, \mu}=0$,
for all $\la, \mu$.

With this fact now established, the same argument as in Proposition
3.1 works well here, if we simply replace the polynomials
$p_{\la,\mu}(z)$ in $z$ with scalars $p_{\la,\mu}$. \QED

\ni 3.3. {\bf{Lemma}}\ \  {\it Let $V$ be a submodule of
$M_{\varphi}$. Then $V$ contains a nonzero Whittaker vector.}

\ni{\bf{Proof}} \ \ Suppose $V$ contains no nonzero Whittaker
vectors. Use the notation as in 2.3.3.

Now let $n_0 = max\{ mindeg(v) | v\neq 0, v \in V \},\, l = min
\{\ell (v) | mindeg (v) = n_0, v \in V \}$.  Take a $u \in V$ such
that $mindeg (u)= n_0$, $\ell (u) = l$. Write $u = u_0w + u_{-1}w
+\cdots + u_{n_0}w, u_i \in \uu''_i, u_{n_0} \neq 0$. Since $u$ is
not a Whittaker vector, there exists a $x \in (\ww_{+})_m$, for some
$m
> 0$ such that $u' := xu-\varphi(x)u = \sum\limits_{i=n_0}^{0} [x, u_i]w \neq 0$. Note $u'$ is contained in
$V$. Clearly, by 2.2.3 Lemma B, $mindeg( [x, u_i]w) \geq i \geq n_0
,$ if $[x, u_i] \neq 0$. So we have  $mindeg(u') \geq n_0$ and hence
$mindeg(u') = n_0 $ for the definition of $n_0$. In this case, $[x,
u_{n_0}]w \neq 0$ and $mindeg([x, u_{n_0}])=n_0$. But this forces $
\ell ([x, u_{n_0}]w ) < ht(u_{n_0}) = \ell(u)$ by Lemma B. Thus,
$\ell ( u' ) < l $, which contradicts with the definition of $l$.
\QED

\vs{18pt}

\cl{\bf\S4. \ Simple Whittaker Modules}

Now we are ready to determine all the simple Whittaker modules.
$\ww, \varphi : \ww_+ \rar \C, M = M_{\varphi}, w = \bo \otimes 1,$
as above.\vspace{10pt}

\ni 4.1 {\bf propostion} \ \ {\it Any nontrivial submodule of a
Whittaker module of type $\varphi$ contains a nontrivial Whittaker
submodule of type $\varphi$.}

\ni{\bf Proof} \ \ Let $V$ be a Whittaker module of type $\varphi$.
We first show it for the case $V = M/IM$ where $I$ is an ideal of
$A=S(z)$. Note that $V = M/IM$ admits a basis, $\{ L_\la
W_\mu\bar{w} \,|\, \la,\mu \in \p- \}$ over $A/I$. Namely, note that
$M = \bigoplus\limits_{\la,\mu \in \p-}AL_\la w$. Hence,
$$M/IM = A/I \otimes_A M = A/I \otimes_A (\bigoplus_{\la,\mu\in
\p-}Ax_\la w)=\bigoplus_{\la, \mu\in \p-}(A/I) x_\la \bar{w} .$$
Then applying the argument in the proof of Lemma 3.3, one sees
immediately that the proposition holds in this case.

Note that with the fact that $V = M/IM$ admits a basis, $\{ L_\la
W_\mu\bar{w} \,|\, \la,\mu \in \p- \}$, over $A/I$, one sees that
Proposition 3.1 holds for $V=M/IM$. Obviously it is enough to show
that there are no other cases. Now, the proposition follows
immediately from the claim below.

 \ni {\bf Claim}: \ \ {\it Let $N$ be a submodule of $M =
M_{\varphi}$. Then $N = IM$ for some $I \subseteq A = S(z)$.}

{ Proof of the claim}: \ \ Set $I = \{ x \in A \,|\, xw \in N \}$.
One immediately sees that $I$ is an ideal of $A$ and $IM \subseteq
N$. So we can view $N/IM$ as a submodule of $M/IM$. If $N\neq IM$,
then there exists $p\bar{w} \in N/IM$, with $p\bar{w} \neq 0, p \in
A$, since 3.1 and 3.3 hold for $M/IM$. So $pw \in N$ and hence $p
\in I$. Therefore $pw \in IM$, which contradicts with the fact that
$p\bar{w} \neq 0$ in $N/IM$. Thus, $N = IM$.\QED\vspace{10pt}

\ni4.2. {\bf Theorem} \ \ {\it For any $\xi \in \C$, $L_{\varphi,
\xi}$ is simple and any simple Whittaker module of type $\varphi$ is
of form $L_{\varphi,\xi}$.}

\ni{\bf Proof} \ \ The fist statement follows from Proposition 3.2
and 4.1. For the second one, let $V$ be a simple Whittake module of
type $\varphi$. Consider a surjection
$$\pi: M \rightarrow V.$$
By Schur's Lemma, there exists a $\xi\in \C$ such that $zv=\xi v,
\forall v\in V$. Hence, $\pi((z-\xi)M_{\varphi})=0$ i.e. $\pi$
factor through $L_{\varphi, \xi}$ and hence $V \simeq L_{\varphi,
\xi}$.\QED\vspace{10pt}

\ni 4.3. We develop two more results to close this
section.\vspace{6pt}

\ni4.3.1. {\bf Lemma} \ \ {\it Set $\uu = \uu (\ww )$,
$L=\uu(z-\xi\cdot1)+\sum\limits_{\la,\mu \in \pq}\uu(L_\la
W_\mu-\varphi(L_\la W_\mu))$, and $V=\uu/L.$ Then $V\simeq
L_{\varphi,\xi}$.}

\ni{\bf Proof} \ \ Note that $\bar{1}$ in $V$ is obviously a
Whittaker vector of type $\varphi$, and also $z$ acts on $V$ by some
scalar $\xi$. By the universal property of $M_{\varphi}$, we have a
surjection
  $$\psi: M_{\varphi}\rightarrow V ,$$
sending $w$ to $\bar{1}$. But then $\psi((z-\xi)M )=(z-\xi)V=0.$
Hence
  $$(z-\xi)M_{\varphi}\subseteq \mbox{ker}\psi \subsetneq M_{\varphi},$$
and therefore, $(z-\xi)M_{\varphi}= \mbox{ker}\varphi$, i.e.
$V\simeq L_{\varphi,\xi}. $ \QED\vspace{6pt}

\ni 4.3.2. {\bf Proposition} \ \ {\it Suppose $V$ is a Whittaker
module, and $z$ acts on $V$ by the scalar $\xi\in \C$. Then $V$ is
isomorphic to $L_{\varphi, \xi}$ and therefore, if $w$ is a cyclic
Whittaker vector for $V$, ${Ann}_{\uu
(\ww)}(w)=\uu(\ww)(z-\xi\cdot1)+\mbox{$\sum\limits_{\la,\mu \in
\pq}$}\uu(\ww )(L_\la W_\mu-\varphi(L_\la W_\mu))$.}

\ni{\bf Proof} \ \ Let $\pi: \uu(\ww)\rightarrow V$, with $\pi
(1)=w$ be the canonical homomorphism, and $K=\mbox{Ker}\pi$. Then
$K\subsetneq \uu(\ww)$ and
$$L:=\uu(\ww)(z-\xi\cdot1)+\mbox{$\sum\limits_{\la, \mu }$}\uu(\ww)(L_\la
W_\la-\varphi(L_\la W_\la))\subset K.$$ By 4.3.1, $L$ is maximal and
thus $K=L$ and $V\simeq \uu(\ww)/L\simeq L_{\varphi, \xi}$. \QED

\vs{18pt}

\cl{\bf\S5. \ Submodules Of Whittaker Modules}

We now characterize arbitrary Whittaker modules, with generating
Whittaker vector $w$, in terms of the annihilator
$\mbox{Ann}_{S(z)}(w)$. The results and their proofs here are
essentially same as \cite{OW}.\vspace{10pt}

\ni 5.1. {\bf \ Decomposition of Whittaker Modules}.\vspace{6pt}

\ni 5.1.1. {\bf Lemma } \ \ {\it Suppose that $V$ is a Whittaker
module of type $\varphi$ with cyclic Whittaker vector $w$ and assume
that $\mbox{Ann}_{S(z)}(w)=(z-\xi\cdot1)^a$ for some $a>0$. define
$V_i\triangleq \uu(\ww)(z-\xi\cdot1)^iw,\ \ 0\leq i\leq a$. Then

1) $V=V_0\supseteq V_1\supseteq\cdots\supseteq V_a=0$ form a
composition series with $V_i/V_{i+1}\simeq L_{\varphi, \xi}$;

2) $V_0\cdots V_a$ are all the submodules of $V$.}

\ni{\bf Proof}\ \ 1) Clearly $V=V_0\supseteq
V_1\supseteq\cdots\supseteq V_a=0$. Since $z$ acts by the scalar
$\xi$ on $V_i/V_{i+1}$, it follows by Proposition 4.3.2 that
$V_i/V_{i+1}$ is simple and isomorphic to $ L_{\varphi, \xi}$.

2) If $M$ is any maximal submodule of $V$, then $V/M$ is simple and
it's easy to use Proposition 4.3.2 to deduce that $V/M \simeq
L_{\varphi, \xi}$. So $(z-\xi\cdot1)V\subset M$, i.e. $V_1\subset
M.$ Therefore $V_1=M$. A similar argument shows that $V_{i+1}$ is
the unique maximal submodule of $V_i, \forall i< a$.\QED\vspace{6pt}

\ni 5.1.2. {\bf Theorem } \ \ {\it Suppose that $V$ is a Whittaker
module of type $\varphi$ with cyclic Whittaker vector $w$, and $
\mbox{Ann}_{S(z)}(w)\neq0$. Let $p(z)$ is the unique monic generator
of $\mbox{Ann}_{S(z)}(w)$. Write
$p(z)=\prod\limits_{i=1}^k(z-\xi_i\cdot1)^{a_i},\ \ \xi_i\neq\xi_j,\
\ i\neq j,$ and $w_j=p_j(z)w $, where $ p_j(z)=\prod\limits_{i\neq
j}(z-\xi_i\cdot1)^{a_i}$. Then $ V_j :=\uu (\ww)w_j$ is a Whittaker
module of type $\varphi$ with cyclic Whittaker vector $w_i$ and
$V=V_1\oplus\cdots\oplus V_k$. Furthermore, the submodules $V_1,
\cdots, V_k$ are indecomposable and the composition length of $V_j$
is $a_j$.}

\ni{\bf Proof} \ \ Since $\mbox{gcd}(p_1(z),\cdots,p_k(z))=1$, there
exist polynomails $q_i(z), 1\leq i \leq k$ such that $\sum
q_i(z)p_i(z)=1$. Therefore $1\cdot w=\sum q_i(z)p_i(z)w\in
V_1+\cdots+ V_k$ and thus, $V=V_1+\cdots+ V_k$. To show that the sum
is direct, first note that $p_j(z)w_i=0, i\neq j$. Hence
$$w_i=1\cdot w_i=q_i(z)p_i(z)w_i.$$
Now if $u_1 w_1+\cdots+u_k w_k=0$, then
$$0=q_i(z)p_i(z)(\mbox{$\sum\limits_j$}u_jw_j)=u_iq_i(z)p_i(z)w_i=u_iw_i.$$
and this implies that the sum is direct. The rest follows from 3.3.3
since $\mbox{Ann}_{s(z)}(w_i) = (z-\xi_i)^{a_i}$.

\ni {\bf Remark}. It's easy to see, from the proof above, that 1)
any submodule of $V$ is the direct sum of its intersections with the
$V_j$'s; 2) $(z-\xi_i)V_j = V_j, i\neq j$.\vspace{6pt}

\ni 5.1.3. {\bf Corollary } \ \ {\it Suppose that $V$ is a Whittaker
module of type $\varphi$ with cyclic Whittaker vector $w$, and $
\mbox{Ann}_{S(z)}(w)=(p(z))$, where $p(z)$ is a monic polynomial.
Then $\mbox{Ann}_{\uu(\ww)}(w)=\uu(\ww)p(z)+\mbox{$\sum\limits_{\la,
\mu }$}\uu(\ww)(L_\la W_\mu-\varphi(L_\la W_\mu))$.}

\ni{\bf Proof}. \ \ By induction on $deg(p(z))$. Assume $deg(p) >1$.
Write $p(z)=(z-\xi\cdot1)p'(z), p'(z)\neq 1$, and $p'(x)$ is monic.
Then $(z-\xi\cdot1)w\neq0$. Let $V'=\uu(\ww)w'$ and
$w'=(z-\xi\cdot1)w$. Then obviously
 $V'$ is a Whittaker module, and $\mbox{Ann}_{S(z)}w'=S(z)p'(z)$. By
induction, we have $\mbox{Ann}_{\uu (\ww )}(w')=\uu (\ww )
p'(z)+\mbox{$\sum\limits_{\la, \mu }$}\uu(\ww)(L_\la
W_\mu-\varphi(L_\la W_\mu)$. Observe that
$\mbox{Ann}_{S(z),V/V'}(\overline{w})=S(z)(z-\xi\cdot1)$, where $
\overline{w}=w+V'$. Since $\mbox{Ann}_{\uu (\ww )}(w)\subset
\mbox{Ann}_{\uu (\ww )}(\overline{w})$, So for any $u\in
\mbox{Ann}_{\uu (\ww )}(w)$, there exist $u_0, u_{\la, \mu} \in \uu
(\ww )$ such that
$$u=u_0(z-\xi\cdot1)+\mbox{$\sum\limits_{\la, \mu\in \pq
}$}u_{\la,\mu}(L_\la W_\mu-\varphi(L_\la W_\mu).$$ But
$\mbox{$\sum\limits_{\la, \mu \in \pq}$} u_{\la, \mu}(L_\la
W_\mu-\varphi(L_\la W_\mu)\in \mbox{Ann}_{\uu (\ww )}(w)$. So
$u_0(z-\xi\cdot1)\in \mbox{Ann}_{\uu (\ww)}(w)$. Thus,
$0=u_0(z-\xi\cdot1)w=u_0w'$, and hence $u_0\in \mbox{Ann}_{\uu (\ww
)}(w')=\uu (\ww )p'(z)+\mbox{$\sum\limits_{\la, \mu }$}\uu (\ww
)(L_\la W_\mu-\varphi(L_\la W_\mu)$. This implies that $u \in \uu
(\ww)p(z)+\mbox{$\sum\limits_{\la, \mu }$}\uu (\ww )(L_\la
W_\mu-\varphi(L_\la W_\mu)$. \QED \vspace{10pt}

\ni 5.2. \ \ {\bf Submodules of $M_{\varphi}$}.

\ni 5.2.1. {\bf Lemma} \ \ {\it Suppose that $V$ is a Whittaker
module of type $\varphi$ with cyclic Whittaker vector $w$. If
$\mbox{Ann}_{S(z)}(w)=0$, then $V\simeq M_\varphi$.}

\ni{\bf Proof} \ \ Consider a surjection, $\pi: M_\varphi\rightarrow
V$. Let $K=\mbox{Ker}\pi$. If $K\neq 0$, then by Lemma 3.3, there
exists $w'\in K$ so that $w'$ is a nonzero Whittaker vector of type
$\varphi$. Hence by 3.1, we have $w'=p(z)( \bo\otimes1)$ for some
nonzero polynomial $p(z)$, and then $0 = \pi(p(z)(\bo
\otimes1))=p(z)w$. Therefore, $p(z)\in \mbox{Ann}_{S(z)}w=0$.
Impossible. So, $K=0$. \QED\vspace{6pt}

\ni 5.2.2 {\bf Theorem} \ \ {\it Set $M = M_\varphi, w=\bo
\otimes1$. Let $V$ be a submodule of $M$. Then $V\simeq M$.
Furthermore, $V$ is generated by a Whittaker vector of the form
$q(z)w$ for some polynomial $q(z)$.}

\ni{\bf Proof} \ \ If $V=M$, then done. Now assume $V \neq M$. Note
first that $\uu (\ww )f(z)w$ is a Whittaker module with cyclic
Whittaker vector $f(z)w$, for any nonzero polynomial $f(z)$. Thus
$\uu (\ww) f(z)w \simeq M$ by Lemma 5.2.1, since
$\mbox{Ann}_{S(z)}(w)=0$. Now Proposition 3.1 and Lemma 3.3 imply
that there exist a nonzero polynomial $p(z)$ such that $w':= p(z)w$
is a Whittaker vector contained in $V$. Therefore
$$M' :=  \uu (\ww )p(z)w\simeq M.$$

We operate on the triple $(M' \subseteq V \subseteq M)$. Write
$p(z)=\mbox{$\prod_{i=1}^r$}(z-\xi_i)^{a_i}, a_i>0, \xi_i\neq\xi_j$.
Denote by $\psi$ the canonical map $M \rightarrow M/M'$, and then
$V=\psi^{-1}(\psi(V))$ since $V\supset M'$. Applying Theorem 5.1.2
to $M/M'$, we have that
$$M/M'=\overline{M_1}\oplus \cdots\oplus\overline{M_r}$$
where $\overline{M_i}=\uu (\ww )\overline{w_i}$,and $\overline{w_i}
= \prod_{j\neq i}(z-\xi_j\cdot1)^{a_j} \overline{w}$, with $
\overline{w}=\psi(w)$. By Remark 5.1.2, $\psi (V) = \oplus_{i} (\psi
(V) \cap \bar{M_i})$. Then, without any loss, we may assume
$\psi(V)\cap \overline{M_1}\neq\overline{M_1}$.  By Lemma 5.1.1,
$(\psi(V)\cap \overline{M_1}) \subseteq
 (z-\xi_1\cdot1)\overline{M_1}$. Let $N =
 (z-\xi_1\cdot1)\overline{M_1}+\overline{M_2}+\cdots+\overline{M_r}$
 and hence $\psi (V) \subseteq N $.  Using Remark 5.1.2, we can
 deduce that $N = (z-\xi_1)M/M'$. Therefore $ M' \subseteq V
 \subseteq \psi^{-1}(N)$. However $\psi^{-1}(N) = (z-\xi_1 )M \simeq
 M$.  With this isomorphism, we obtain another triple $(M'' \subseteq
 V' \subseteq M )$ with $ M'' \simeq M' , V' \simeq V$.
 Observe that the composition length of $M/M''$  is one
 less than the one of $M/M'$. Therefore, to complete the proof, we only
  need to repeat the above operation on the new triple. \QED

Remark.\, One may also use the facts developed in the proof of
Proposition 4.1 to give another proof.

\ni5.3.\ \ Now we can characterize  submodules of any Whittaker
modules.

\ni 5.3.1. {\bf Lemma } \ \ {\it Suppose that $V$ is an
indecomposable Whittaker module of type $\varphi$ with cyclic
Whittaker vector $w$, then $Wh(V)=S(z)w$. where $Wh(V)$ stands for
the set of all Whittaker vectors (including $0$) of V.}

\ni{\bf Proof} \ \ a) $\mbox{Ann}_{S(z)}(w)=0$, the result follows
from Lemma 5.2.1 and Proposition 3.1.

b) $\mbox{Ann}_{S(z)}(w)\neq0$. Use induction on $\ell (V)$, the
composition length of $V$. By 5.1, we have
$$V=V_0\supseteq V_1\supseteq\cdots\supseteq V_a=0.$$
where $V_i$ is generated by the cyclic Whittaker vector
$w_i=(z-\xi\cdot1)^iw$ for some $\xi \in \C$, and $a = \ell (V)$.
For $w'\in Wh(V)$,  we conclude, by 5.1.1 and 3.2, that
$\overline{w'}=c\overline{w}$ in $V/V_1$ for some $c\in \C$.
Therefore $w'=cw+w''$, with $w''\in V_1$. Note $w''=w'-cw$ is also a
Whittaker vector. Now $\ell (V_1) = a-1 $ and then by induction,
$w''=p(z)w_1=p(z)(z-\xi\cdot1)w$ for some $p(z)\in S(z)$, therefore
$w'=cw+p(z)(z-\xi\cdot1)w$. So $Wh(V) = S(z)w$.\QED

\ni5.3.2 {\bf Corollary} \ \ {\it Suppose that $V$ is a  Whittaker
module of type $\varphi$. Then 1) there exist a series of
indecomposable Whittaker modules $V_i, 1\leq i \leq r $, such that
$V = \oplus_{i=1}^{r}V_i$; 2) Suppose $w$ is a cyclic Whittaker
vector of $V$. Then any submodule of $V$ is a Whittaker module
generated by a vector in $Wh(V) = S(z)w$.}

\ni {\bf Proof} \ \ 1) follows from Theorem 5.1.2 and Lemma 5.2.1.
2) If $V \simeq M_{\varphi}$,then it follows from Theorem 5.2.2; if
not, it follows from 5.1 and 5.3.1. \QED

\end{document}